\documentclass[a4paper,11pt]{amsart}

\usepackage{amssymb,amscd,amsmath}
\usepackage[mathcal]{eucal}
\usepackage{array,float}

\usepackage{pdflscape}
\usepackage{hyperref}
\usepackage{pdflscape}
\usepackage{enumerate}
\usepackage{relsize}
\usepackage{xcolor}

\usepackage[left=3cm, right=3cm]{geometry}
\usepackage{enumerate}
\numberwithin{equation}{section}
\numberwithin{equation}{subsection}

\newtheorem{theorem}{Theorem}[section]

\newtheorem{corollary}[theorem]{Corollary}
\newtheorem{lemma}[theorem]{Lemma}

\newtheorem*{remark*}{Remark}
\newtheorem{remark}[theorem]{Remark}

\title[Finite groups with mostly involuted cyclic subgroups]{Finite groups with mostly involuted cyclic subgroups} 

\author{Vaibhav Chhajer}

\address{School of Mathematical Sciences, National Institute of Science Education and Research Bhubaneswar, 752050, Odisha, India}
\address{Homi Bhabha National Institute, Training School Complex, Anushakti Nagar, Mumbai, 400094, India}

\email{vchhajer@niser.ac.in}
\email{vaibhavchhajer0504@gmail.com }

\author{Palash Sharma}
\address{Department of Mathematical Sciences, Indian Institute of Science Education and Research Mohali, Knowledge City, Sector 81, Mohali 140 306, Punjab, India}

\email{ms22001@iisermohali.ac.in}
\email{palashsharma8765434124@gmail.com}

\begin{document}
\date{}
\begin{abstract} 
Let $G$ be a finite group. Let $C(G)$ be the set of cyclic subgroups of $G$, $c(G)=|C(G)|$ and $i(G)=\left|\{x\in G : x^{2}=e\}\right|$. In this article, we classify finite groups with $i(G)=c(G)-r$ for $r\in \{0,1,2\}$ and prove that the range of the function given by $\beta(G)=\frac{i(G)}{c(G)}$ is dense in $[0,1]$. We also answer an open question posed by Gao and Shen in `Finite groups with many cyclic subgroups'.

\end{abstract}
\maketitle

\section{Introduction}

Throughout this article, $\mathbb{Z}_{n}$, $D_{2n}$, $S_n$, and $A_n$ denote the cyclic group of order $n$, the dihedral group of order $2n$, the symmetric group, and the alternating group on $n$ letters, respectively. Suppose $\text{Aut}(G)$ denotes the automorphism group of a group $G$, $\phi$ denotes the Euler totient function, and $\tau$ denotes the divisor function. Let $\exp(G)$ denote the exponent of a group $G$. Let $|G|$ and $o(x)$ denote the order of a group $G$ and the order of an element $x \in G$, respectively. Let $G_{1}\rtimes G_{2}$ denote a semidirect product of a group $G_{1}$ with another group $G_{2}$. Let $n_{p}$ denote the number of Sylow $p$-subgroups of a group $G$. Let $S(G)$ be the set of subgroups of $G$ and $s(G)=|S(G)|$. Let $\mathcal{F}$ denote the set of all finite groups.

Classifying finite groups based on the number of involutions has always been an interesting problem in group theory. In \cite{M}, G.A. Miller studied groups in which at least half of their elements were involutions. Then, in \cite{EN}, the authors fully classified finite groups with at least three-fourths of their elements as involutions. 

Recently, the problem of classifying finite groups based on the number of cyclic subgroups present in them has also been considered. In \cite{MT}, the author asked the following question and answered it for $r=1$:

\begin{quote}
    \textbf{Open problem:} Describe the finite groups $G$ with $c(G) = |G|-r$, where $2 \le r \le |G|-1$.
\end{quote}
This problem has been studied by several authors in \cite{AA, AA1, RJL, LCQ,CHS, SR1, SR, MT1}, and motivated by these results, researchers also began studying the function $\alpha(G) = \frac{c(G)}{|G|}$ in \cite{GS, GL, ML}. Following a similar idea, and combining the concepts of classifying finite groups based on the number of involutions and the number of cyclic subgroups, we characterize finite groups $G$ in which most of the cyclic subgroups are of order at most $2$. Our first main result classifies finite groups by analyzing the equation $i(G) = c(G) - r$ for small values of $r$, and is stated as follows:

\begin{theorem}\label{mainresult}
    Let $G$ be a group with $i(G)=c(G)-r$ and $p$ be an odd prime, then
    \begin{enumerate}
        \item $r=0$ if and only if $G$ is isomorphic to $\mathbb{Z}_{2}^{n}$, for $n\geq 0$.
        \item $r=1$ if and only if $G$ is isomorphic to one of the following groups: $\mathbb{Z}_4, D_8, \mathbb{Z}_p$, or $D_{2p}$.
        \item $r=2$ if and only if $G$ is isomorphic to one of the following groups: $\mathbb{Z}_2\times\mathbb{Z}_{4}, \mathbb{Z}_2\times D_{8}, \mathbb{Z}_8, D_{16}, \mathbb{Z}_{p^2}, D_{2p^2}, \mathbb{Z}_{2p}$, or $D_{4p}$.
    \end{enumerate}
\end{theorem}

The subgroup lattice of a finite group $G$ and its associated posets contain structural information and help to understand various properties of a group (see \cite{RS1}). Recalling the definition of an atom from order theory, we observe that within the poset of cyclic subgroups $C(G)$, the subgroups of order $2$ are specific atoms. In this context, we observe in Lemma \ref{lemma4.4} that if $\alpha(G) \ge \frac{3}{4}$, meaning at least $\frac{3}{4}$ of the elements in a finite group $G$ generate distinct cyclic subgroups, then at least $\frac{2}{3}$ of the cyclic subgroups of $G$ are of order at most $2$.

To further analyze the relationship between the number of involutions and the number of cyclic subgroups, it is natural to consider the function $\beta: \mathcal{F}\to [0,1]$ given by $\beta(G)=\frac{i(G)}{c(G)}$, which measures the probability of a random cyclic subgroup of $G$ being of order at most $2$. From a lattice-theoretic standpoint, the invariant $\beta(G)$ precisely quantifies the cardinality of the subposet formed by the trivial subgroup and the involution-generated atoms relative to the entire poset $C(G)$.

Similar to $\alpha$ and $\beta$, there are several other group-theoretic functions, such as $\gamma:\mathcal{F}\to [0,\infty)$ and $\operatorname{cdeg}: \mathcal{F}\to [0,1]$, defined as  $\gamma(G)=\frac{s(G)}{|G|}$ and $\operatorname{cdeg}(G)=\frac{c(G)}{s(G)}$. Determining the density of the ranges of such functions has always been a topic of interest. For instance, the author in \cite[Theorem 4.3]{SML} showed that the range of $\gamma$ is dense in  $[0,\infty)$, and the same author later showed that the range of $\operatorname{cdeg}$ is dense in $[0,1]$ in \cite[Theorem 3.2]{RS}. In contrast, the range of $\alpha$ is not dense in $[0,1]$, as is evident from \cite[Theorem 5]{GL}. On similar lines, our second main result establishes the density of the range of $\beta$:

\begin{theorem}\label{mainresult2}
    The set $S=\{\beta(G): G\in \mathcal{F} \}$ (the range of $\beta$) is dense in $[0,1]$.
\end{theorem}

Consequently, the density of the range of $\beta$ in $[0, 1]$ implies that the subposet formed by subgroups of order at most $2$ can be realized in any proportion of the entire poset $C(G)$. This demonstrates that the combinatorial structure of the cyclic subgroup poset imposes no rigid restrictions on the relative size of this subposet.

Additionally, in Remark \ref{remark}, we answer an open question regarding group solvability posed in \cite{GS}.

\section{Prerequisites}

 First, we recall the following results of group theory, which will be used to prove Theorem \ref{mainresult}.
\begin{lemma}\label{normal}
    Let $G$ be a finite group and $H\leq G$. If $H$ is a unique cyclic subgroup of order $|H|$, then $H\trianglelefteq G$. 
\end{lemma}

\begin{theorem}[\text{\cite[Theorem 4.1]{EN}}]\label{3/4 inv}
 If $G$ is a finite group and $i(G)> \frac{3}{4}|G|$, then $G$ is an elementary
 abelian 2-group.
\end{theorem}

\begin{theorem}\label{c(G)=|G|-r} 
    Let $G$ be a finite group with $c(G)=  |G|-r$. Then
    \begin{enumerate}
        \item $r=1$  if and only if $G$  is isomorphic to one of the following groups:  $\mathbb{Z}_3, \mathbb{Z}_4, S_3,$ or $D_8$ \text{(see  \cite[Theorem 2] {MT})}.  

        \item  $r=2$ if and only if  $G$  is isomorphic to one of the following groups:  $\mathbb{Z}_{6}, \mathbb{Z}_{4}\times \mathbb{Z}_2, D_{12},$ or $ \mathbb{Z}_2\times D_{8}$ \text{(see \cite[Theorem 2]{RJL})}. 

        \item  $r=4$ if and only if $G$  is isomorphic to one of the following groups:  $\mathbb{Z}_{4}\times \mathbb{Z}_{2}\times \mathbb{Z}_{2}, \mathbb{Z}_{2}\times \mathbb{Z}_{2}\times D_8, (\mathbb{Z}_{2}\times \mathbb{Z}_{2})\rtimes \mathbb{Z}_{4}, Q_{8}\rtimes\mathbb{Z}_{2}, \mathbb{Z}_{3}\times\mathbb{Z}_{3}, (\mathbb{Z}_{3}\times\mathbb{Z}_{3})\rtimes \mathbb{Z}_{2}, A_{4}, \mathbb{Z}_{6}\times\mathbb{Z}_{2}, \mathbb{Z}_{2}\times\mathbb{Z}_{2}\times S_3, \mathbb{Z}_{8}, $ or $ D_{16}$ \text{(see \cite[Theorem 4]{RJL})}. 
    \end{enumerate}
   
\end{theorem} 

 \begin{theorem}\label{p^k or 2p^k}
     Let $p$ be an odd prime and $n\in\{p^{k},2p^{k}\}$ for some $k\geq 1$. Then $\mathbb{Z}_{n}\rtimes \mathbb{Z}_{2}$  is isomorphic to one of the following groups:  $\mathbb{Z}_{n}\times\mathbb{Z}_{2}$ or $D_{2n}$.
 \end{theorem}
 \begin{proof}
     The semidirect product $G = \mathbb{Z}_n \rtimes \mathbb{Z}_2$ induces a group homomorphism $\varphi : \mathbb{Z}_2 \to \operatorname{Aut}(\mathbb{Z}_n)$. If $\varphi$ is the trivial homomorphism, then the semidirect product is simply the direct product, and we are done since $G \cong \mathbb{Z}_n \times \mathbb{Z}_2$. If $\varphi$ is non-trivial, let $h$ denote the generator of $\mathbb{Z}_2$. It follows that $\varphi(h)$ must be an automorphism of order $2$ of $\mathbb{Z}_n$. Hence, $\varphi(h)(g) = g^t$ for all $g \in \mathbb{Z}_n$, where $t^2 \equiv 1 \pmod{n}$ and $t \neq 1$. Recall that $\operatorname{Aut}(\mathbb{Z}_m)$ is cyclic if and only if $m \in \{1, 2, 4, q^k, 2q^k\}$ where $k \ge 1$ and $q$ is an odd prime. Therefore, $\operatorname{Aut}(\mathbb{Z}_n)$ is a cyclic group of even order, so it has a unique element of order $2$. So, there is only one non-trivial semidirect product. Without loss of generality, we may assume that $t = -1$. This leads to $hgh^{-1} = g^{-1}$ for all $g \in \mathbb{Z}_n$. In particular, if $g$ generates $\mathbb{Z}_n$, then $G$ is generated by $g$ and $h$, and the following equalities hold: 
\[
g^n = e, \quad h^2 = e, \quad \text{and} \quad hgh^{-1} = g^{-1}.
\]

So it follows that $G \cong D_{2n}$.
 \end{proof}

We also recall the following results of analysis, which will be used to prove Theorem \ref{mainresult2}
\begin{theorem}[\text{\cite[Theorem 1.1]{NSH}}]\label{limit test} 
    Let $\sum_{n=1}^\infty a_n$ and $\sum_{n=1}^{\infty} b_n$ be series with positive terms. If $$\lim_{n\to \infty} \frac{a_{n}}{b_{n}}=c$$ where $c$ is a finite positive number, then either both series are convergent or both are divergent.
\end{theorem}

\begin{lemma}[\text{\cite[Lemma 4.1]{SML}}]\label{denseness result} 
     Let $(x_n)_{n\geq 1}$ be a sequence of positive real numbers such that $\lim_{n\to \infty} x_{n}=0$ and the series $\sum _{i=1}^{\infty} x_{i} $ is divergent. Then, the set containing the sums of all finite subsequences of $(x_{n})_{n\geq 1}$ is dense in $[0,\infty).$
\end{lemma}

\section{Proof of the first main result}
In this section, we prove our first main result.
But before that, we present the following basic result:
\begin{lemma}\label{basic results 1}
    (a) Let $n$ be a positive integer. Then $\mathbb{Z}_n$ and $D_{2n}$ both satisfy $i(G)=c(G)-r$, where 
    \begin{equation*}
        r=
        \begin{cases} 
       \tau(n)-1 & \text{if } n \text{ is odd,} \\
      \tau(n)-2 & \text{if } n \text{ is even.}
   \end{cases}
    \end{equation*}
    
    (b) If a group $H$ is a solution of $i(G)=c(G)-r$, then $H\times \mathbb{Z}_2$ is a solution of $i(G)=c(G)-2r.$ 
\end{lemma}
\begin{proof}
    The proof of part $(a)$ is standard and therefore omitted. For part $(b)$, suppose $H$ is a solution to the equation $i(G) = c(G) - r$, meaning $i(H) = c(H) - r$. Utilizing the properties of the direct product with $\mathbb{Z}_2$, namely $i(H \times \mathbb{Z}_2) = 2i(H)$ and $c(H \times \mathbb{Z}_2) = 2c(H)$, we obtain:
\[
i(H \times \mathbb{Z}_2) = 2i(H) = 2(c(H) - r) = 2c(H) - 2r = c(H \times \mathbb{Z}_2) - 2r.
\]
Therefore, $H \times \mathbb{Z}_2$ is a solution to $i(G) = c(G) - 2r$.
\end{proof}

\textbf{Proof of Theorem \ref{mainresult}.}
\begin{proof} Set $c=c(G)$ and $i=i(G)$. Clearly, $i\leq c$. Without loss of generality, let $x_1,x_2,\cdots,x_i,x_{i+1},\cdots,x_c$ be the generators of cyclic subgroups of $G$, where $x_1=e$ and $x_2,\cdots,x_i$ are involutions.

For $r=0$, it can be easily seen that this occurs if and only if $G \cong \mathbb{Z}_2^n$, where $n \geq 0$.

When $r>0$, we observe that $G\ncong \mathbb{Z}_2^n$ and $i=|G|-(\phi(o(x_{c-(r-1)}))+\cdots+\phi(o(x_{c})))$. Therefore, by Theorem \ref{3/4 inv} we get $$|G|\leq 4(\phi(o(x_{c-(r-1)}))+\cdots+\phi(o(x_{c}))).$$ Also $c=|G|-(\phi(o(x_{c-(r-1)}))+\cdots+\phi(o(x_{c}))-r)$.

For $r=1$, we have  $c=|G|-(\phi(o(x_{c}))-1)$. Now  using the fact that for each divisor $k$ of $o(x_{c})$, there is a  cyclic subgroup of  order $k$ sitting inside $\langle x_{c}\rangle$, we conclude that $o(x_{c})\in\{ p,4\}$, where $p$ is an odd prime.

If $o(x_{c})=4$, then $c=|G|-1 $ and so by Theorem \ref{c(G)=|G|-r}, $\mathbb{Z}_{4}$ and $ D_{8}$ are solutions.

If $o(x_{c})=p$, then  $|G|=2^{k}p$ for some $k\geq 0$ and $|G|\leq 4(p-1)$. Therefore $|G|\in\{p,2p\}$. 
\begin{enumerate}
    \item If $|G|=p$, then $G\cong \mathbb{Z}_{p}$ is a solution.

    \item If $|G|=2p$, then by Lemma \ref{normal}, $\langle x_{c}\rangle \trianglelefteq G$ and so by Theorem \ref{p^k or 2p^k}, we get $D_{2p}$ as solution.
\end{enumerate}

For $r=2$, we have $c=|G|-(\phi(o(x_{c-1}))+\phi(o(x_{c}))-2)$. From \cite[Section 3.3.3]{CHS} we have the following possibilities:

\begin{center}
\begin{tabular}{ |c|c|c|c|c|c|c| } 
\hline
&(a) & (b) & (c) & (d) &  (e) & (f) \\
\hline
$o(x_{c-1})$ & $4$ & $4$  & $p$ & $p$ &  $p$& $p$\\ 
$o(x_{c})$ & $4$ & $8$ & $p^2$ & $2p$  & $4 \hspace{1mm}\text{or}\hspace{1mm} q$ & $p$\\ 
\hline
\end{tabular}
\end{center} 
where $p$ and $q$ are distinct odd primes. We'll deal with them one by one.

\textbf{(a) Assume $o(x_{c-1})=4$ and $o(x_{c})=4$.} Then $c=|G|-2$  and so by Theorem \ref{c(G)=|G|-r}, $\mathbb{Z}_{4}\times \mathbb{Z}_2$ and $\mathbb{Z}_2\times D_{8}$ are solutions.

\textbf{(b) Assume $o(x_{c-1})=4$ and $o(x_{c})=8$.} Then $c=|G|-4$ and so by Theorem \ref{c(G)=|G|-r}, $\mathbb{Z}_{8}$ and $ D_{16}$ are solutions.

\textbf{(c) Assume $o(x_{c-1})=p$ and $o(x_{c})=p^2$.} Then $|G|=2^{k}p^2$ for some $k\geq 0$ and $ |G|\leq 4(p^2-1)$, therefore $|G|\in\{p^2,2p^2\}$.
\begin{enumerate}
    \item  If $|G|=p^2$, then $o(x_{c})=p^2$ implies that $G\cong \mathbb{Z}_{p^2}$ is a solution.

    \item If $|G|=2p^2$, then by Lemma \ref{normal}, $\langle x_{c}\rangle \trianglelefteq G$ and so by Theorem \ref{p^k or 2p^k}, we get $D_{2p^2}$ as solution.
\end{enumerate}

\textbf{(d) Assume $o(x_{c-1})=p$ and $o(x_{c})=2p$.} Then $|G|=2^{k}p$ for some $k\geq 1$ and  $|G|\leq 8(p-1)$, therefore $|G|\in\{2p,4p\}$.
\begin{enumerate}
\item If $|G|=2p$, then  $o(x_{c})=2p$ implies that $G\cong \mathbb{Z}_{2p}$ is a solution.

    \item If $|G|=4p$, then by Lemma \ref{normal}, $\langle x_{c}\rangle \trianglelefteq G$ and so by Theorem \ref{p^k or 2p^k}, we get $D_{4p}$ as solution. 
\end{enumerate}

\textbf{(e) Assume $o(x_{c-1})=p$ and $o(x_{c})\in\{4 ,q\}$.} Then by Lemma \ref{normal}, $\langle x_{c-1}\rangle, \langle x_{c}\rangle$ $\trianglelefteq G $. Since $\langle x_{c-1}\rangle\cap \langle x_{c}\rangle = \{e\}$, $G$ contains a subgroup isomorphic to $\mathbb{Z}_{4p}$ or $\mathbb{Z}_{pq}$, a contradiction.

\textbf{(f) Assume $o(x_{c-1})=p$ and $o(x_{c})=p$.} Then $|G|=2^{k}p$ for some $k\geq 0$ and $|G|\leq 8(p-1)$, therefore $|G|\in \{p,2p,4p\}$. The cases $|G|\in\{p,2p\}$ are not possible because no group of these orders contains two different cyclic subgroups of order $p$. Hence, $|G|=4p$. If $p>3$, then $n_p\mid 4$ and $n_p \equiv 1 \hspace{1mm} \text{mod} \hspace{1mm}p$ implies that $G$ has a unique Sylow $p$- subgroup, and so $\langle x_{c-1}\rangle = \langle x_{c}\rangle$, a contradiction. Therefore $p=3$ and $|G|=12$, but there are no groups of order $12$ having $2$ distinct subgroups of order $3.$
\end{proof}

The following is an immediate consequence of Theorem \ref{mainresult}.
\begin{corollary}
    If $G$ is a finite group such that $i(G)=c(G)-r$ with $r\in\{0,1,2\}$, then $G$ is a supersolvable group. 
\end{corollary}

We ask a natural open problem related to the above study. 
\vspace{5mm}

\noindent\textbf{Open problem:} 
Classify all finite groups $G$ satisfying $i(G)=c(G)-r$, where $2< r < c(G)$.

\section{Proof of the second main result.}

First, observe the following necessary and sufficient condition for $\beta$ to be multiplicative:

\begin{lemma} \label{basic property of beta}
    Let $G$ and $H$ be finite groups. Then 
    $$\beta(G\times H)=\beta(G)\cdot\beta(H)\ 
    \text{if and only if } \operatorname{gcd}(\operatorname{exp}(G),\operatorname{exp}(H)) \in \{1,2\}.$$
\end{lemma}

\begin{proof}
    Clearly, $i(G\times H) = i(G)\cdot i(H)$ holds unconditionally. Therefore, $\beta(G\times H)=\beta(G)\cdot\beta(H)$ if and only if $c(G\times H)=c(G)\cdot c(H)$. 
    
    Recall that for any finite group $K$, $c(K) = \sum_{x \in K} \frac{1}{\phi(o(x))}$. Thus, 
    $$c(G\times H) = \sum_{(x,y)\in G\times H} \frac{1}{\phi(\operatorname{lcm}(o(x),o(y)))}.$$
    Using the property of $\phi$ that $\phi(\operatorname{lcm}(m,n))\cdot \phi(\operatorname{gcd}(m,n)) =\phi(m)\cdot\phi(n)$ for all $m,n \in \mathbb{N}$, we obtain
    $$c(G\times H) = \sum_{(x,y)\in G\times H} \frac{\phi(\operatorname{gcd}(o(x),o(y)))}{\phi(o(x))\cdot\phi(o(y))}.$$
    Since $\phi(k) \ge 1$ for all $k \in \mathbb{N}$, we always have $c(G\times H) \ge c(G)\cdot c(H)$, with equality holding if and only if $\phi(\operatorname{gcd}(o(x),o(y))) = 1$ for all $x \in G$ and $y \in H$, that is, if and only if $ \operatorname{gcd}(o(x),o(y)) \in \{1,2\}$ for all $x \in G$ and $y \in H$. Therefore, it suffices to show that $ \operatorname{gcd}(o(x),o(y)) \in \{1,2\}$ for all $x \in G$ and $y \in H$ holds if and only if $\operatorname{gcd}(\operatorname{exp}(G),\operatorname{exp}(H)) \in \{1,2\}$. 
    
    If $\operatorname{gcd}(\operatorname{exp}(G),\operatorname{exp}(H)) \in \{1,2\}$, then since $o(x)$ divides $\operatorname{exp}(G)$ and $o(y)$ divides $\operatorname{exp}(H)$, it immediately follows that $\operatorname{gcd}(o(x),o(y))$ divides $\operatorname{gcd}(\operatorname{exp}(G),\operatorname{exp}(H))$, and thus $\operatorname{gcd}(o(x),o(y))$ $\in \{1,2\}$.
    
    Conversely, suppose that $\operatorname{gcd}(o(x), o(y)) \in \{1, 2\}$ for all $x \in G$ and $y \in H$, but, for the sake of contradiction, that $\operatorname{gcd}(\operatorname{exp}(G), \operatorname{exp}(H)) \notin \{1, 2\}$. Then we have the following two cases:
    
    \textbf{Case 1:} Suppose $\operatorname{gcd}(\operatorname{exp}(G), \operatorname{exp}(H)) = 2^n r$, where $n \geq 0$ and $r \geq 3$ is an odd integer.
Then, an odd prime $p$ divides $\operatorname{gcd}(\operatorname{exp}(G), \operatorname{exp}(H))$ and so $p \mid |G|$ and $p \mid |H|$. Thus, there exist $a \in G$ and $b \in H$ such that $o(a) = p$ and $o(b) = p$. Consequently, $\operatorname{gcd}(o(a), o(b)) = p \notin \{1, 2\}$, a contradiction.

  \textbf{Case 2:} Suppose $\operatorname{gcd}(\operatorname{exp}(G), \operatorname{exp}(H)) = 2^n$, where $n \geq 2$. Then, $4 \mid \operatorname{exp}(G)$ and $4 \mid \operatorname{exp}(H)$. 
If it were true that $o(x) \leq 2$ for all $x \in G$, this would imply $\operatorname{exp}(G) \leq 2$, which contradicts $4 \mid \operatorname{exp}(G)$. Therefore, there exists $a \in G$ such that $o(a)=4$. Similarly, there exists $b \in H$ such that $o(b)=4$.
Consequently, $\operatorname{gcd}(o(a), o(b)) = 4 \notin \{1, 2\}$, a contradiction.

    Since both cases lead to a contradiction, we conclude that $\operatorname{gcd}(\operatorname{exp}(G), \operatorname{exp}(H)) \in \{1, 2\}$.
\end{proof}

The following are immediate observations that can be made using the above lemma.

\begin{corollary}\label{coprime or Z2n}
    If $G$ and $H$ are groups with $\operatorname{gcd}(|G|,|H|)=1$ or $H\cong \mathbb{Z}_2^n$, where $n\geq 0$, then 
    $$\beta(G\times H)=\beta(G)\cdot \beta(H).$$
\end{corollary}

\begin{proof}
    If $\operatorname{gcd}(|G|,|H|)=1$, then $\operatorname{gcd}(\operatorname{exp}(G),\operatorname{exp}(H)) = 1$. If $H\cong \mathbb{Z}_2^n$, then $\operatorname{exp}(H) \in \{1,2\}$. Consequently, in either case, $\operatorname{gcd}(\operatorname{exp}(G),\operatorname{exp}(H)) \in \{1,2\}$, and the result follows trivially from Lemma \ref{basic property of beta}.
\end{proof}

\begin{corollary}\label{D_2p and D_2q}
    Let $\{p_j\}_{j=1}^{n}$ be $n$ distinct odd primes. Then 
    $$\beta\left(\prod_{j=1}^{n}D_{2p_j}\right)=\prod_{j=1}^{n}\beta(D_{2p_j}).$$
\end{corollary}

\begin{proof}
    We proceed by induction on $n$. The base case $n=1$ is trivial. For $n \geq 2$, let $G = \prod_{j=1}^{n-1}D_{2p_j}$ and $H = D_{2p_n}$. Then $\operatorname{exp}(G)=2p_1 \cdots p_{n-1}$ and $\operatorname{exp}(H)=2p_n$. Since $p_n$ is a distinct odd prime from $p_1, \dots, p_{n-1}$, we have $\operatorname{gcd}(\operatorname{exp}(G), \operatorname{exp}(H)) = 2$. By Lemma \ref{basic property of beta}, $\beta(G \times H) = \beta(G)\cdot\beta(H)$, and the result follows by the inductive hypothesis.
\end{proof}

\begin{lemma}
    If $G $ is a non-trivial group of odd order, then $ \beta(G) \le \frac{1}{2},$ with equality if and only if $G \cong \mathbb{Z}_p $ for some odd prime $p$.
\end{lemma}
\begin{proof}
Let $G$ be a non-trivial group of odd order. Then $i(G) = 1$. Consequently, we have $\beta(G) = \frac{1}{c(G)}$. Since $G$ is non-trivial, $c(G) \ge 2$. Therefore, $\beta(G) \le \frac{1}{2}$. Equality holds if and only if $c(G) = 2$, and trivially $c(G) = 2$ if and only if $G \cong \mathbb{Z}_p$ for some odd prime $p$.
\end{proof}

\begin{lemma}\label{lemma4.4}
    Let $G$ be a finite group. If $\alpha(G)\ge \frac{3}{4}$, then $\beta(G)\ge \frac{2}{3}$. Moreover, if $\alpha(G)\ge \frac{5}{6}$, then $\beta(G)\in\{ \frac{4}{5},\frac{6}{7},1\}$.
\end{lemma}
\begin{proof}
    According to \cite[2.2]{GL}, we have the inequality $\beta(G) \ge 2 - \frac{1}{\alpha(G)}$. This immediately gives that if $\alpha(G) \ge \frac{3}{4}$, then $\beta(G) \ge \frac{2}{3}$. Furthermore, if $\alpha(G) = \frac{3}{4}$, \cite[Main Theorem]{GS} shows that this lower bound of $\frac{2}{3}$ is exactly achieved by a specific class of groups. Finally, if $\alpha(G) \ge \frac{5}{6}$, then by \cite[Corollary 2]{GL}, we deduce that $\beta(G) \in \left\{ \frac{4}{5}, \frac{6}{7}, 1 \right\}$.
\end{proof}

\begin{remark}\label{remark}
    In \cite{GS}, the authors asked the following question:
    \begin{quote}
    Let $A$ be a finite abelian group and $G$ a finite group. Suppose that $\alpha(G) = \alpha(A)$ and $|G| = |A|$. Is $G$ also solvable?
\end{quote}
\end{remark}

We now answer this question in the negative by providing a counterexample. Let $G = A_5 \times \mathbb{Z}_{16}$ and $A = \mathbb{Z}_{15} \times \mathbb{Z}_2^2 \times \mathbb{Z}_4^2$. Then, $|G| = |A| = 960$ and $\alpha(G) = \alpha(A) = \frac{1}{6}$. However, $G$ is not a solvable group because it contains the alternating group $A_5$ as a subgroup, which is nonsolvable.\\

Now we prove our second main result.\\

\textbf{Proof of Theorem \ref{mainresult2}.}
\begin{proof}
We'll follow the same idea as presented in the proof of \cite[Theorem 3.2]{RS} and consider $\mathbb{R}$ with the standard topology. Let $n$ be a positive integer and $p_n$ denote the $n^{\text{th}}$ odd prime number. Then we have
$$ \beta(D_{2p_n}) = \frac{p_n+1}{p_n+2}. $$
Consider the sequence $(x_n)_{n \ge 1} \subseteq (0, \infty)$, where $x_n = \ln\left(\frac{p_n+2}{p_n+1}\right)$. We have
$$ \lim_{n\to\infty} \frac{x_n}{1/p_n} = \lim_{n\to\infty} \frac{\ln\left(\frac{p_n+2}{p_n+1}\right)}{1/p_n} = 1. $$
As the series $\sum_{n=1}^\infty \frac{1}{p_n}$ diverges, by Theorem \ref{limit test}, we deduce that the series $\sum_{n=1}^\infty x_n$ also diverges. As $\lim_{n\to\infty} x_n = 0$ so by Lemma \ref{denseness result} we have
$$ \overline{\left\{ \sum_{j \in I} \ln\left(\frac{p_j+2}{p_j+1}\right) : I \subseteq \mathbb{N}, |I| < \infty \right\}} = [0, \infty) \iff $$
$$ \overline{\left\{ \ln\left(\prod_{j \in I} \frac{p_j+2}{p_j+1}\right) : I \subseteq \mathbb{N}, |I| < \infty \right\}} = [0, \infty). $$ 

Let $B= {\left\{ \ln\left(\prod_{j \in I} \frac{p_j+2}{p_j+1}\right) : I \subseteq \mathbb{N}, |I| < \infty \right\}}$. Since the function
$ f: \mathbb{R} \to \mathbb{R},$ given by $f(x) = e^{-x}, \forall \ x \in \mathbb{R},$
is continuous, we have:

$$ f(\overline{B}) \subseteq \overline{f(B)} \implies f\left([0,\infty)\right) = (0,1] \subseteq \overline{f(B)} \implies [0,1]\subseteq \overline{f(B)},$$ because $\overline{f(B)}$ is closed in $\mathbb{R}$.

And as $f(B)\subseteq [0,1]\implies \overline{f(B)}\subseteq\overline {[0,1]}=[0,1]$. Hence $\overline{f(B)}=[0,1].$ Therefore, $$\overline{\left\{ \prod_{j \in I} \frac{p_j+1}{p_j+2} : I \subseteq \mathbb{N}, |I| < \infty \right\}}=  \overline{\left\{ \prod_{j \in I} \beta(D_{2p_j}) : I \subseteq \mathbb{N}, |I| < \infty \right\}}= [0,1].$$ Then by using Corollary \ref{D_2p and D_2q} we have $$ \overline{\left\{ \beta\left(\prod_{j\in I}D_{2p_j}\right) : I \subseteq \mathbb{N}, |I| < \infty \right\}}= [0,1].$$ 
As $$\left\{ \beta\left(\prod_{j\in I}D_{2p_j}\right) : I \subseteq \mathbb{N}, |I| < \infty \right\}\subseteq S \subseteq (0,1],$$ we get that $S$ is dense in $[0,1]$ by taking closures of the sets outlined in the previous inclusion.
\end{proof}

The proof of Theorem \ref{mainresult2} immediately implies the following result.

    \begin{corollary}
    Let $\mathcal{M}$ be the set of all finite direct products of distinct dihedral groups of the form $D_{2p}$ where $p$ is an odd prime, and let $\mathcal{C}$ be a set of finite groups such that $\mathcal{M} \subseteq \mathcal{C}$. Then the set $S_\mathcal{C} = \{\beta(G) : G \in \mathcal{C}\}$ is dense in $[0,1]$.
\end{corollary}

Finally, we pose an open question motivated by the preceding result.

\vspace{5mm}

\noindent\textbf{Open problem:} Let $\mathcal{F}$, $C(G),$ and $c(G)$ be as before. Let $n>2$ be an integer. For $G\in \mathcal{F}$, define $i_{n}(G)=\left|\{H\in C(G): |H|\in \{ 1,n\}\}\right|$. Consider the function $\beta_{n}: \mathcal{F}\to [0,1]$ given by $\beta_{n}(G)=\frac{i_n(G)}{c(G)}$,  which measures
the probability of a random cyclic subgroup of $G$ being of order $1$ or $n$. Let $S_{n}=\{\beta_{n}(G): G\in \mathcal{F}\}$ (range of $\beta_{n}$). Is $S_{n}$ dense in $[0,1]$?

\vspace{5mm}

\noindent\textbf{Acknowledgements:} 
The authors are grateful to the anonymous referee for carefully reading the manuscript and providing valuable comments and suggestions. The first-named author acknowledges the research support of the SERB-SRG, Department of Science and Technology (Grant no. SRG/2023/000894), Govt. of India, and thanks the National Institute of Science Education and Research, Bhubaneswar, for providing an excellent research facility. The second-named author thanks the Indian Institute of Science Education and Research, Mohali, for providing an excellent research facility. The authors also thank Dr. Sumana Hatui and Dr. Darshan Nasit for helpful discussions and valuable input.\\

\noindent\textbf{Disclosure statement:}  We have no conflicts of interest to disclose. \\

\noindent\textbf{Author contributions:} Both authors contributed equally to this work.


\begin{thebibliography}{999}







\bibitem{AA}
A.R. Ashrafi, and E. Haghi, On the number of cyclic subgroups in a finite group, Southeast Asian Bull. Math., 42 (2018), 865-873.

\bibitem{AA1}
A.R. Ashrafi, and E. Haghi, On n-cyclic groups, Bull. Malays. Math. Sci. Soc., 42 (2019), 3233-3246.

\bibitem{RJL}
R. Belshoff, J. Dillstrom, and L. Reid, Finite groups with a prescribed number of cyclic subgroups, Commun. Algebra, 47.3 (2019), 1043-1056.

\bibitem{LCQ}
X. Chen, H. Liu, and S. Qiao, Finite groups with a small number of cyclic subgroups, Czechoslovak Math. J., 75(150) (2025), 839-851.

\bibitem{CHS}
V. Chhajer, S. Hatui, and P. Sharma, Finite groups with nearly half as many cyclic subgroups as elements, arXiv:2506.21163v2.

\bibitem{EN}
A. L. Edmonds, and Z. B. Norwood, Finite groups with many involutions, arXiv:0911.1154.

\bibitem{GS}
X. Gao, and R. Shen, Finite groups with many cyclic subgroups, Indian J. Pure Appl. Math., 54(2) (2023), 485-498.

\bibitem{GL}
M. Garonzi, and I. Lima, On the number of cyclic subgroups of a finite group, Bull. Braz. Math. Soc. (N.S.), 49(3) (2018), 515-530.

\bibitem{NSH}
N. S. Hoang, A limit comparison test for general series, Am. Math. Mon., 122.9 (2015), 893-896.

\bibitem{SML}
S. M. Lazorec, A connection between the number of subgroups and the order of a finite group, J. Algebra Appl., 21(1) (2022), Paper No. 2250001.

\bibitem{RS}
S. M. Lazorec, Element orders in extraspecial groups, Acta Math. Hungar., 173.2 (2024), 434-447.

\bibitem{ML}  
S. M. Lazorec, and M. Tărnăuceanu, A note on the number of cyclic subgroups of a finite group, Bull. Math. Soc. Sci. Math. Roumanie (N.S.), 62(110) (2019), 403-416.

\bibitem{M} 
G. A. Miller, Groups Containing a Relatively Large Number of Operators of Order Two, Bull. Am. Math. Soc., 25(9) (1919), 408-413.

\bibitem{RS1}
R. Schmidt, Subgroup Lattices of Groups, De Gruyter Expositions in Mathematics 14, De Gruyter, Berlin (1994).

\bibitem{SR1}
K. Sharma, and A.S. Reddy, Groups having $11 $ cyclic subgroups, Int. J. Group Theory, 13(2) (2024), 203-214.

\bibitem{SR}
K. Sharma, and A.S. Reddy, Groups having $12 $ cyclic subgroups, Bull. Aust. Math. Soc., 113(2) (2026), 304-315.

\bibitem{MT}
M. Tărnăuceanu, Finite groups with a certain number of cyclic subgroups, Am. Math. Mon., 122.3 (2015), 275-276.

\bibitem{MT1}
M. Tărnăuceanu, Finite Groups with a Certain Number of Cyclic Subgroups II, Acta Univ. Sapientiae Math., 10(2) (2018), 375-377.

\end{thebibliography}
\end{document}